\documentclass[12pt,twoside]{article}
\usepackage{amssymb,amsmath,amstext,amsthm,amsfonts}
\usepackage[ansinew]{inputenc} 
\usepackage{graphicx}
\usepackage[mathscr]{eucal}
\usepackage[pdfnewwindow=true]{hyperref}

\newcommand{\N}{\mathbb{N}}

\newcommand{\Mod}[1]{\left\vert{#1}\right\vert}

\newcommand{\prl}[2]{\langle\!\langle {#1};{#2}\rangle\!\rangle }
\newcommand{\pair}[2]{[\![\overline{#1},\overline{#2}]\!] }
\newcommand{\lgpair}[2]{\left[\!\left[{#1},{#2}\right]\!\right] }
\newcommand{\maxord}[1]{\mbox{maxord}({#1})}
\newcommand{\ord}{\mbox{ord}}

\newcommand{\coment}[1]{}

\newcommand{\Ascr}{\mathscr{A}}
\newcommand{\Dscr}{\mathscr{D}}

\newtheorem{defi}{Definition}
\newtheorem{lema}{Lema}
\newtheorem{prop}{Proposition}
\newtheorem{teor}{Theorem}

\newcommand{\dem}{ \par\medbreak\noindent{\bf
Proof. }\enspace} 
 
\newcommand{\cqd}{\hfill
$\sqcup\!\!\!\!\sqcap\bigskip$}

\title{A discrete Faà di Bruno's formula}
\author{ P. Duarte$\,^{1}$,  \; M. J. Torres$\,^{2}$ }
\date{}

\begin{document}
\maketitle

{\small
\noindent
$^1$ {\it CMAF, Departamento de Matem\'{a}tica, Faculdade de Ci\^{e}ncias da Universidade de Lisboa, Portugal
(e-mail: pedromiguel.duarte@gmail.com).
Supported by Funda\c{c}\~{a}o para a Ci\^{e}ncia e a Tecnologia, Financiamento Base 2008 - ISFL/1/209.} \\
$^2$ {\it CMAT, Departamento de Matem\'{a}tica, Universidade do Minho, Campus de Gualtar, 4700-057 Braga,
Portugal (e-mail: jtorres@math.uminho.pt). }
}

\bigskip

\noindent
{\bf Abstract }
{\em 
We derive some formulas that
rule the behaviour of finite differences under
composition of functions with vector values and arguments.
}

\bigskip

\section{Introduction}

The Faà di Bruno formula~\cite{Fa} gives an expression for the $n$-th derivative ($ n\geq 1$)
of the composition\, $f\circ g$ \, of two functions $f$ and $g$
in terms of derivatives of $f$ and $g$.
Functions of class $C^\infty$ will be called smooth. 
To simplify assumptions we shall always assume that $f$ and $g$ are both smooth functions. 
This formula has many versions, depending on the type of derivatives considered.
Assume first that $f$ and $g$ are real valued functions of one real argument.
Then the formula takes the following form

\begin{equation} \label{faadibruno:1d}
(f\circ g)^{(n)} (x) =
 \sum  \; \frac{n!}{b_1!\, \cdots\, b_n!}\,f^{(b_1+\cdots +b_n)}(g(x))\, 
\left( \frac{g'(x)}{1!} \right)^{b_1}\,\cdots\,
\left( \frac{g^{(n)}(x)}{n!} \right)^{b_n}
\end{equation}
where the sum is taken over all solutions $(b_1,\ldots, b_n)\in \N^n$ 
of the equation 
$$b_1+2\,b_2+\ldots + n\, b_n = n\;.$$
We note that $\N=\{0,1,\ldots\}$ starts with zero, so that $b_i\geq 0$ for every $i=1,\ldots, n$.

Assume next that $f$ and $g$ are functions between euclidean spaces of possibly different
dimensions, the domain of $f$ containing the range of $g$. Let $e_1,\ldots, e_k$ be vectors
in the euclidean domain of $g$. The list $e=(e_1,\ldots, e_k)$ will be referred as a
$k$-multi-vector. Given $\alpha=(\alpha_1,\ldots, \alpha_k) \in \N^k$ and a point $x$ in the
domain of $g$, we denote by $D^\alpha g_x(e)$ the derivative
$D^{\alpha_1}_{e_1}\circ \cdots \circ D^{\alpha_k}_{e_k}  g(x)$
of order $\Mod{\alpha}=\alpha_1+\ldots+\alpha_k$.
Given $\alpha\in\N^k$, let 
$$[\alpha]=\{\, \beta\in\N^k\,:\, \beta\leq \alpha\,\} =
\{ 0, \beta^1, \beta^2,\ldots, \beta^m\}\;.$$
Then, Faà di Bruno's formula takes the form
\begin{equation} \label{faadibruno:md}
D^\alpha  (f\circ g)_x(e)  = 
 \sum  \; \frac{\alpha!}{b_1!\, \cdots\, b_m!}\,
D^{(b_1,\ldots ,b_m)} f_{g(x)}\, 
\left( \frac{D^{\beta^1} g_x(e) }{(\beta^1)!}, \, \ldots\,,\,
\frac{D^{\beta^m} g_x(e)}{(\beta^m)!} \right) 
\end{equation}
where the sum is taken over all solutions $(b_1,\ldots, b_m)\in \N^m$ 
of the equation 
$$b_1 \beta^1 + b_2 \beta^2+\ldots + b_m \beta^m = \alpha\;,$$
and $\alpha!= \alpha_1 !\,\alpha_2 ! \,\cdots\,\alpha_k!$ when $\alpha=(\alpha_1,\ldots, \alpha_k)$.
Assume now the multi-index
$\alpha$ belongs to the discrete cube $I^k=\{0,1\}^k$.
In this case we have  $\alpha!=1$, \, $\beta^i!=1$,\, and $b_i\in \{0,1\}$,
which implies $b_i !=1$, for all $i=1,\ldots, m$.
Then formula~(\ref{faadibruno:md}) reduces to

\begin{equation*}  
D^\alpha  (f\circ g)_x(e)  = 
 \sum  \; 
D^{(b_1,\ldots ,b_m)} f_{g(x)}\, 
\left( D^{\beta^1} g_x(e) , \, \ldots\,,\,
D^{\beta^m} g_x(e) \right) 
\end{equation*}
where the sum is taken over all solutions $(b_1,\ldots, b_m)\in \N^m$ 
of the equation 
$$b_1 \beta^1 + b_2 \beta^2+\ldots + b_m \beta^m = \alpha\;.$$
We can simplify this formula a little more.
\begin{defi}\label{partition}
Given a multi-index $\alpha\in I^k$, we  call
{\em partition of $\alpha$ with size $r$} to any subset
with $r$ elements $\{\alpha^1\,\ldots,\alpha^r\}\subseteq I^k-\{0\}$
such that $\alpha^1+\cdots +\alpha^r = \alpha $.
\end{defi}
Of course the size $r$ must range between $1$ and $\Mod{\alpha}$, since
the smallest partition $\{\alpha\}$ has size $1$, and $\Mod{\alpha^i}\geq 1$
for each $\alpha^i$ in some partition. Then

\begin{equation} \label{faadibruno:I}
D^\alpha  (f\circ g)_x(e)  = 
 \sum_{\alpha^1 +\cdots + \alpha^r = \alpha}   \; 
D^{m} f_{g(x)}\, 
\left( D^{\alpha^1} g_x(e) , \, \ldots\,,\,
D^{\alpha^r} g_x(e) \right) 
\end{equation}
where the sum is taken over all partitions of $\alpha$ with size
ranging from $1$ to $\Mod{\alpha}$.
Formula~(\ref{faadibruno:I}) is at the same time a special case 
and an extension of~(\ref{faadibruno:md}).
Given $\alpha\in\N^k$, take $\tilde{e}$ to be the
$r$-multi-vector with  $\alpha_i$ components equal to $e_i$,
for each $i=1,\ldots, k$, where $r=\Mod{\alpha}$.
Set $\tilde{\alpha}\in I^r=\{0,1\}^r$ to be the multi-index with all components equal to $1$.
Then \, $D^\alpha  (f\circ g)_x(e) =
D^{\tilde{\alpha}}  (f\circ g)_x(\tilde{e}) $. Grouping and counting equal terms
we can derive~(\ref{faadibruno:md}) from~(\ref{faadibruno:I}).

\bigskip

After having written this article we found the reference~\cite{FL}
where a discrete version of this formula is derived in terms of divided differences.
Here we deal with {\em finite differences} in stead of {\em divided differences},
which means our formula holds for general functions of several variables.
The history Finite Difference Calculus goes back a long way,
parallel to that of Infinitesimal Calculus. 
We refer to~\cite{J} for a modern treatment of this calculus.
We briefely recall some basic definitions in order to
state our discrete version of Faà di Bruno's formula.
Given a vector $u\in X$, in some euclidean space $X$, let $\tau_u$ be the 
{\em translation operator} defined by $(\tau_u f)(x)=f(x+u)$.
This operator acts on every space of functions $f\in Y^X$.
The {\em difference operator along vector $u$},\, $\Delta_u:Y^X\to Y^X$, is defined by\,
$\Delta_u={\rm id} -\tau_u$, which corresponds to set\,
 $(\Delta_u f)(x)=f(x+u)-f(x)$.
Notice these operators always commute, i.e., 
$\Delta_u\circ \Delta_v = \Delta_v\circ \Delta_u$\,
for all vectors $u, v\in X$, since \,
$\tau_v\circ\tau_u
 = \tau_{u+v } = \tau_u\circ\tau_v$.
\begin{defi}\label{finite:diff}
We call finite difference operator of order $k$ to any composition of $k$
difference operators along possibly repeated vectors.
Given a multi-vector $u=(u_1,\cdots,u_k)\in X^k$, 
we denote by
$\Delta^k_u :Y^X \rightarrow Y^X$ the composite operator
$\Delta^k_u = \Delta_{u_1}\circ\Delta_{u_2}\circ\cdots\circ \Delta_{u_k}$.
More generally, given $\alpha\in \N^k$,
$\Delta^\alpha_u = \left(\Delta_{u_1}\right)^{\alpha_1} \circ \ldots \circ 
\left(\Delta_{u_k} \right)^{\alpha_k} $
denotes a difference operator of order $\Mod{\alpha}$.
\end{defi}

Next we introduce an algebra of symbolic finite difference expressions.
Let $x,y,\ldots$ be symbols representing points,
$u_1, u_2, \ldots$ be symbols representing vectors,
and $f, g, \ldots$ be symbols representing functions.
We denote by $\Dscr$ the set of all symbolic finite difference expressions,
which we define as the smallest set of expressions such that:
\begin{enumerate}
\item if $x$ is a symbol representing a point then $x\in\Dscr$;
\item if $u$ is a symbol representing a vector then $u\in\Dscr$;
\item if $t \in\Dscr$  and $f$ is a symbol representing a function then
$f(t) \in\Dscr$;
\item if $t, s\in\Dscr$ then  $t+s \in\Dscr$;
\item if $\alpha\in \N^k$, $s, t_1,\ldots, t_k \in\Dscr$  and $f$ is a symbol representing a function then
$\Delta^{\alpha}_{(t_1,\ldots, t_k)} f(s) \in\Dscr$.
\end{enumerate}
We consider as equal all terms which formally can be proved to be equal using 
property transformation rules of finite differences.
Of course, depending on the interpretation given to the point, vector and function symbols,
many terms in $\Dscr$ will be formal but meaningless expressions.
We define recursively the order of a term ${\rm ord}:\Dscr\to \N$:
\begin{enumerate}
\item if $x$ is a symbol representing a point then~ ${\rm ord}(x)=0$;
\item if $u$ is a symbol representing a vector then~ ${\rm ord}(u)=1$;
\item if $t \in\Dscr$  and $f$ is a symbol representing a function then~
${\rm ord}(f(t))=0$;
\item if $t, s\in\Dscr$ then~  ${\rm ord}(t+s)=\min\{{\rm ord}(t), {\rm ord}(s) \} $;
\item if $\alpha\in \N^k$, $s, t_1,\ldots, t_k \in\Dscr$  and $f$ is a symbol representing a function then 
${\rm ord}\left(\Delta^{\alpha}_{(t_1,\ldots, t_k)} f(s)\right)=
\alpha_1\,{\rm ord}(t_1)+\ldots + \alpha_k\,{\rm ord}(t_k)$.
\end{enumerate}
Given a term $t=\Delta^{\alpha}_{(t_1,\ldots, t_k)} f(s)\in\Dscr$ which is meaningful
for some interpretation of its symbols (all functions being smooth),
if all vectors $u_i$ in $t$ are small of order \, $\epsilon$\, then\,
$t$ is small of order $\epsilon^{{\ord}(t)}$.
The following theorem is our main result.
Let $X$, $Y$ and $Z$ stand for linear spaces.

\bigskip


\noindent
{\bf Theorem A} {\em 
Given  maps $f\in Y^X$ and $g\in X^Z$, and a
multi-vector  $u\in Z^k$,
\begin{equation}
\label{dcomp}
 \Delta^\alpha_u (f\circ g)(x) = \sum_{\alpha^1+\cdots+\alpha^n=\alpha}
\Delta^{\,n}_{\Delta^{\alpha^1}_u g(x), \ldots , \Delta^{\alpha^n}_u g(x) } 
f\left(\,g(x)\, \right)\; +\;\cdots\; \;,
\end{equation}
where the ellipsis stand for higher order terms.

}

\bigskip

Next theorem gives an explicit formula for~(\ref{dcomp}).

\bigskip

\noindent
{\bf Theorem B} {\em 
Given $\alpha\in I^k$, there are recursively defined sets, 
$\Ascr^\xi_0, \, \Ascr^\xi_{\alpha^1},\,\ldots,\, \Ascr^\xi_{\alpha^r}$
associated with each
partition 
$\xi=\{\alpha^1,\,\ldots,\,\alpha^r\}\subseteq I^k$ of $\alpha$, such that
\begin{enumerate}
\item the sets $\Ascr^\xi_0, \, \Ascr^\xi_{\alpha^1},\,\ldots,\, \Ascr^\xi_{\alpha^r}$ are pairwise disjoint;
\item $\beta\in \Ascr^\xi_{\beta}\subset [\alpha]$, for $\beta= 0, \alpha^1,\ldots, \alpha^r$;
\item $\Mod{\gamma}>\Mod{\beta}$, for every $\gamma\in \Ascr^\xi_\beta-\{\beta\}$; and 
\end{enumerate}
for any given  maps $f\in Y^X$ and $g\in X^Z$, and any
multi-vector  $u\in Z^k$,
\begin{equation}
\label{tda}
 \Delta^\alpha_u (f\circ g)(x) = \sum_{\xi = \{\alpha^1,\ldots,\alpha^r\}\in {\mathscr P}_\alpha}
\Delta^{\,r}_{ u^\xi_{\alpha^1},\,\ldots,\, u^\xi_{\alpha^r} } f(u^{\xi}_0) \;,
\end{equation}
where
$u^{\xi}_\beta = \sum_{\gamma\in \Ascr^\xi_\beta} \Delta^\gamma_u g(x)$,\,
for each $\beta= 0, \alpha^1,\ldots, \alpha^r$.
}

\bigskip

See theorem~\ref{TeorTdiscrfexpr}.
Note that 
$u^{\xi}_0 = g(x) + \cdots$, with remainder $\sum_{\gamma\in \Ascr^\xi_0-\{0\}} \Delta^\gamma_u g(x)$,
and for each $i$,\,
$u^{\xi}_{\alpha^i} = \Delta^{\alpha^i}_u g(x) + \cdots $,\, with remainder\,
$\sum_{\gamma\in \Ascr^\xi_{\alpha^i}-\{\alpha^i\}} \Delta^\gamma_u g(x)$.
By 3. both these remainders are terms of  higher order.

\bigskip

\section{The Infinitesimal Formula}

Let $D\subseteq  X$ be an open domain. The tangent space $T(D)$, 
and the tangent map $Tf:T(D)\to T(Y)$,
are defined to be
$T(D)=D\!\times\! X$,\, respectively  \, $Tf(x,u)=\left(f(x), D_u f(x)\right)$.
The chain rule shows that construct $T$ is a functor, which essentially means that
$T(f\circ g)=Tf \circ Tg$. Inductively, we can define higher order tangent
spaces and tangent maps, by setting
$T^k(D)=T(T^{k-1}(D))$ and $T^k(f)=T(T^{k-1}(f))$. 
Then the iterated correspondence $T^k$ becomes also a functor.
The tangent map of order $k$, $T^k(f)$ can be explicitly expressed
in terms of higher order directional derivatives of $f$.
As we shall see the pattern of these expressions rules the behaviour
of higher order derivatives under composition. 
To get an explicit expression for $T^k(f)$ we need some
special notation to denote elements in $T^k(D)$.
We shall call {\em $k$-cuboid} of $X$ to any family
$\overline{u}=(u_\alpha)_{\alpha\in I^k}$ 
of vectors in $X$ indexed over the discrete cube $I^k=\{0,1\}^k$.
Notice that any $k$-cuboid $\overline{u}$
can be thought of as a pair of $(k-1)$-cuboids, corresponding
to restrict its indices to the two opposite faces
$\{\alpha_k=0\}$ and $\{\alpha_k=1\}$ of the discrete
cube $I^k$. Therefore, we can and shall  identify the tangent space
$T^k(X)$ with the set of all $k$-cuboids of $X$:
$$ T^k(X)=\{ \, \overline{u}=(u_\alpha)_{\alpha\in I^k}\; :\; u_\alpha\in X\;
\text{ for all }\; \alpha\in I^k\; \}\;.$$
The $k$-tangent space  $T^k(D)$, 
to an open domain $D$, is the set of all $k$-cuboids
$\overline{u}\in T^k(X)$ whose base point $u_0$ belongs to $D$.
The $k$-tangent space over a point $x\in D$ is the
set $T^k_x(D)$ of all $\overline{u}\in T^k(D)$ such that
$u_0=x$.
We shall use the multi-index derivative notation
$$ D^\alpha_u f (x) = \left(D_{u_1}\right)^{\alpha_1} \circ \ldots \circ \left(D_{u_k} \right)^{\alpha_k} f (x)\;. $$

Because multi-indices are cumbersome to write, we shall adopt the following
writing convention. Given a cuboid $\overline{u}\in T^k(X)$,
$u_{i_1,\ldots,i_n}$ stands for the component $u_\alpha$, where
$\alpha$ is the multi-index
$(\alpha_1,\ldots,\alpha_k)$ defined by
$\alpha_j=1$, if $j\in\{i_1,\ldots,i_n\}$, $\alpha_j=0$ otherwise.
Given a multi-index $u=(u_1,\ldots,u_k)\in X^k$
we shall denote by\,
$\prl{x}{u}$\, the $k$-cuboid $\overline{w}$ such that
$w_{0}=x$, $w_{1}=u_1$, $\ldots$, 
$w_{k}=u_k$,
and $w_{i_1,\ldots, i_n}=0$ for all $n\geq 2$.
With this notation is very easy to check that

\begin{prop}
\label{Tfprlexpr}
If $f\in Y^X$ is a class $C^k$ map, then
for all $x\in X$, $u\in X^k$,
$$ T^k(f)\left(\prl{x}{u} \right)= \left( D^\alpha_u f(x) \right)_{\alpha\in I^k }\;.$$
\end{prop}

\bigskip

The following kind of notation
$$ \sum_{\alpha^1+\cdots+\alpha^n=\alpha} A_{\alpha^1,\ldots, \alpha^n} $$
will always denote a sum taken over all partitions of $\alpha$ (see definition~\ref{partition})
with size ranging from $1$ to $\Mod{\alpha}$.

\begin{prop}
\label{Tfexpr}
Given a class $C^k$ map $f\in Y^X$ and a $k$-cuboid $\overline{u}\in T^k_x(X)$,
writing 
$T^k f\, \overline{u} = ( T^\alpha f \, \overline{u} )_{\alpha\in I^k}$,
we have for each $\alpha\in I^k$
\begin{equation}
\label{tde}
 T^\alpha f \, \overline{u} = \sum_{\alpha^1+\cdots+\alpha^n=\alpha}
D^{\,n}_{u_{\alpha^1}, \ldots , u_{\alpha^n} } f(x) \;.
\end{equation}
\end{prop}

\dem
This proposition is proved by induction.
\cqd

Let us call {\em order} of a term
$D^{\,n}_{u_{\alpha^1}, \ldots , u_{\alpha^n} } f(x)$
to the sum
$\Mod{\alpha^1}+\ldots +\Mod{\alpha^n}$.
Then, it is clear that expression~(\ref{tde}) is
homogeneous: all summands have order $\Mod{\alpha}$.

A couple of examples
\begin{eqnarray*}
 T^{(1,1)}(f)(\overline{u})
&=& D_{u_{1,2}} f(x) + D^2_{u_{1},\,u_{2}} f(x) 
\end{eqnarray*}
\begin{eqnarray*}
 T^{(1,1,1)}(f)(\overline{u})
&=& D_{u_{1,2,3}} f(x) + \\
& & D^2_{u_{1},\,u_{2,3}} f(x) + D^2_{u_{2},\,u_{1,3}} f(x) + D^2_{u_{3},\,u_{1,2}} f(x) +\\
& & D^3_{u_{1},\,u_{2},\,u_{3}} f(x). 
\end{eqnarray*}

From proposition~\ref{Tfprlexpr}, we see that
the tangent map pattern~(\ref{tde}) rules the behaviour
of higher order derivatives under composition. 

\begin{prop}
\label{invder}
Given class $C^k$ maps $f\in Y^X$ and $g\in X^Z$, and a
multi-vector  $u\in Z^k$,
\begin{equation*}
 D^\alpha_u (f\circ g)(x) = \sum_{\alpha^1+\cdots+\alpha^n=\alpha}
D^{\,n}_{D^{\alpha^1}_u g(x), \ldots , D^{\alpha^n}_u g(x) } f\left(g(x)\right) \;.
\end{equation*}
\end{prop}

\dem
\begin{eqnarray*}
D^\alpha_u (f\circ g)(x) &=& T^\alpha(f\circ g) \prl{x}{u} \\
&=& T^\alpha(f)\, T^k(g) \prl{x}{u} \\
&=& T^\alpha(f)\, \left( D^\beta_u g(x)\right)_{\beta\in I^k} \\
&=& \sum_{\alpha^1+\cdots+\alpha^n=\alpha}
D^{\,n}_{D^{\alpha^1}_u g(x), \ldots , D^{\alpha^n}_u g(x) } f(x)
\end{eqnarray*}
\cqd

\bigskip

\section{The Discrete Formula}

In order to characterize finite difference operators consider as before the 
discrete cube  $I^k=\{0,1\}^k$ as 
a set of multi-indices. Given $\alpha\in I^k$ 
we write $\alpha=(\alpha_1,\cdots,\alpha_k)$, where each
$\alpha_i$ represents a binary digit, $\alpha_i=0$ or $\alpha_i=1$.
The set $I^k$ is partially ordered by the relation
$$\alpha\leq \beta \quad \Leftrightarrow \quad \alpha_i\leq \beta_i,\;
\text{ for all }\; i=1,\ldots,k\;.$$
We also write
$\Mod{\alpha}=\alpha_1+\ldots +\alpha_k$\; and \;
$\alpha\cdot u = \alpha_1\,u_1+\ldots + \alpha_k\,u_k$.
A simple computation shows that

\begin{prop}
\label{finite:diff:expr}
Given $f\in Y^X$, for all $x\in X$, and $u\in X^k$,
$$\Delta^\alpha_u f (x)= \sum_{\beta\leq\alpha} (-1)^{|\alpha|-|\beta|} f\left(x+\beta\cdot u\right) $$
\end{prop}

\bigskip

A key property of the difference operator $\Delta$ is the following
kind of additivity.
\begin{prop}
\label{diff:add}
Given $f\in Y^X$, for all $x\in X$, and $u,v\in X$,
\begin{equation}
\label{aditiv}
 \Delta_{u+v} f(x) = \Delta_u f(x) + \Delta_v f(x+u)\;.
\end{equation}
\end{prop}
\bigskip

We define now a discrete equivalent of
the $k$-tangent map to $f:X\to Y$.
This will be a mapping $T_k(f):T^k(X)\to T^k(Y)$.
The construct $T_k$ will again be a functor.
For that purpose we introduce the {\em difference operator}
$\Delta:T^k(X)\rightarrow T^k(X)$ 
$$\Delta \overline{u} =( \Delta^\alpha \overline{u} )_{\alpha\in I^k}\;,
\quad\text{ where }
\quad \Delta^\alpha \overline{u}= \sum_{\beta\leq\alpha} (-1)^{|\alpha|-|\beta|} u_\beta\;.$$
Notice that $\Delta^0 \overline{u} = u_0$. 
The operator $\Delta$ is invertible.
Its inverse is the {\em sum operator}
$\Delta^{-1}:T^k(X)\rightarrow T^k(X)$ 
$$ \Delta^{-1}\overline{u} =\left( \, \sum_{\beta\leq\alpha} u_\beta \,\right)_{\alpha\in I^k}\;.$$
The correspondence $f\in Y^X\,\leadsto\, f_\ast : T^k(X)\to T^k (Y)$,
$$f_\ast (x_\alpha)_{\alpha\in I^k} =\left(f(x_\alpha)\right)_{\alpha\in I^k}\;,$$
is obviously a functor.
Thus, defining $T_k (f) = \Delta\circ f_\ast \circ \Delta^{-1}$,
the  correspondence $f\leadsto T_k(f)$
is conjugated to $f\leadsto f_\ast$. Therefore,
$T_k$ behaves functorially too.

\bigskip

\begin{prop}
\label{Tdiscrfprlexpr}
Given $f\in Y^X$,
for all $x\in X$, and $u\in X^k$,
$$ T_k(f)\left(\prl{x}{u} \right)= \left( \Delta^\alpha_u f(x) \right)_{\alpha\in I^k }\;.$$
\end{prop}

\dem
It is enough to notice that
$\Delta^{-1}\prl{x}{u} = (x+\alpha\cdot u)_{\alpha\in I^k}$.
\cqd

\bigskip

We shall say that any component $u_\alpha$, of a $k$-cuboid
$\overline{u}$, has order $\Mod{\alpha}$.
Then, we define recursively the {\em order of a  finite difference 
term } $\Delta^{\,n}_{u_{1}, \ldots , u_{n} } f(x)$
to be  the sum of the orders of terms $u_1$, $\ldots$, $u_n$.
Each term $u_i$ can either be some cuboid component, 
as in proposition below, or else another finite difference term,
as in theorem A.
In both cases the formulas for $T_\alpha f \overline{u}$ and
$\Delta^\alpha_u(f\circ g)(x)$, respectively,
have a main part which is a sum of order $\Mod{\alpha}$ terms,
plus a remainder consisting of terms with order $>\Mod{\alpha}$.

\begin{prop}
\label{Tdiscrfexpr}
Given $f\in Y^X$ and $\overline{u}\in T^k_x(X)$,
writing 
$T_k f\, \overline{u} = ( T_\alpha f \, \overline{u} )_{\alpha\in I^k}$,
we have for each $\alpha\in I^k$
\begin{equation}
\label{tdc}
 T_\alpha f \, \overline{u} = \sum_{\alpha^1+\cdots+\alpha^n=\alpha}
\Delta^{\,n}_{u_{\alpha^1}, \ldots , u_{\alpha^n} } f(\,x \,) \; +\; \cdots \;
\end{equation}
where the dots stand for a sum of higher order terms.
\end{prop}

\dem
This proposition follows from
theorem~\ref{TeorTdiscrfexpr},
using the property~(\ref{aditiv})
to expand differences.
\cqd

\bigskip

Given  $\alpha\in I^k$, we shall denote by
${\mathscr P}_\alpha$ the set of all partitions
$\xi=\{\alpha^1,\ldots,\alpha^n\}$ of $\alpha$.
The multi-index obtained from $\alpha$ adding
digit $1$ at the end will be denoted by $\alpha\Diamond 1$.
Therefore $\alpha\Diamond 1$ belongs to $I^{k+1}$
and has order $\Mod{\alpha\Diamond 1}=\Mod{\alpha}+1$.
Next lemma relates ${\mathscr P}_\alpha$ with
${\mathscr P}_{\alpha\Diamond 1}$.

\begin{lema}
\label{diamond}
Given $\xi=\{\alpha^1,\,\ldots,\,\alpha^n\}\in {\mathscr P}_\alpha$, 
consider the partitions of $\alpha\Diamond 1$
\begin{tabbing}
\quad \=$\tilde{\xi}_0 = \{ 0\Diamond 1,\,\alpha^1\Diamond 0,\,\ldots,\,\alpha^n\Diamond 0\}$
\; with size $n+1$,\\
	\>$\tilde{\xi}_i = \{ \alpha^1\Diamond 0,\,\ldots,\,\alpha^i\Diamond 1,\,\ldots, \alpha^n\Diamond 0\}$
\; with size $n$, \; for\; $1\leq i\leq n$.
\end{tabbing}
These partitions exhaust  ${\mathscr P}_{\alpha\Diamond 1}$ without repetitions,
$$ {\mathscr P}_{\alpha\Diamond 1}=\{ 
\tilde{\xi}_i\, :\; 0\leq i\leq n,\;\text{ and }\; \xi\in {\mathscr P}_\alpha\; \}\;.$$ 
\end{lema} 

\bigskip

Given two $k$-cuboids $\overline{u},\,\overline{v}\in T^k(X)$,
we shall denote by $\overline{w}=\pair{u}{v}$ the unique 
$(k+1)$-cuboid such that
$\overline{w}_{\alpha\Diamond 0} = u_\alpha $ and  
$\overline{w}_{\alpha\Diamond 1} = v_\alpha $.

\begin{lema} For all $\overline{u},\,\overline{v}\in T^k(X)$,
\begin{enumerate}
\item $\Delta^{-1}\pair{u}{v}=\lgpair{\Delta^{-1} \overline{u}}{ \,\Delta^{-1} (\overline{u}+ \overline{v})  }$.
\item $\Delta\pair{u}{v}=\lgpair{\Delta \overline{u}}{ \,\Delta \overline{v}- \Delta \overline{u}  }$.
\item $T_{k+1}(f)\pair{u}{v}=\lgpair{T_k(f) \overline{u}}{ \,T_k(f)(\overline{u} +\overline{v})- T_k(f) \overline{u}  }$.
\end{enumerate}
\end{lema}

\bigskip

\begin{lema}
\label{lemsplit}
\begin{equation*}
\begin{split}
&\Delta^n_{(u_1+v_1),\,\ldots,\,(u_n+v_n)} f(x+w) 
-\Delta^n_{u_1,\,\ldots,\,u_n} f(x) = \\
&\quad  \Delta^{n+1}_{w,u_1,\,\ldots,\,u_n} f(x) +
\sum_{i=1}^n \Delta^n_{u_1,\,\ldots,\,u_{i-1},\,v_i,\,(u_{i+1}+v_{i+1}),\,
\ldots,\,(u_n+v_n)} f(x+w+u_i)
\end{split}
\end{equation*}
\end{lema}

\dem
Follows by an iterated application of~(\ref{aditiv}).
\cqd

\bigskip

We call {\em maximum order} of a partition 
$\xi=\{\alpha^1,\,\ldots,\,\alpha^n\}\in{\mathscr P}_\alpha$
to the number \; $\maxord{\xi}=\max\{\Mod{\alpha^1},\,\ldots,\,\Mod{\alpha^n}\}$.

\bigskip

\begin{teor}
\label{TeorTdiscrfexpr}
Given $\alpha\in I^k$, there are recursively defined sets, 
$\Ascr^\xi_0, \, \Ascr^\xi_{\alpha^1},\,\ldots,\, \Ascr^\xi_{\alpha^r}$
associated with each
partition 
$\xi=\{\alpha^1,\,\ldots,\,\alpha^r\}\subseteq I^k$ of $\alpha$, such that
\begin{enumerate}
\item the sets $\Ascr^\xi_0, \, \Ascr^\xi_{\alpha^1},\,\ldots,\, \Ascr^\xi_{\alpha^r}$ are pairwise disjoint;
\item $\beta\in \Ascr^\xi_{\beta}\subset [\alpha]$, for $\beta= 0, \alpha^1,\ldots, \alpha^r$;
\item $0 <\gamma <\alpha$ \, and \, $\Mod{\gamma}< \maxord{\xi}$, \,
for every $\gamma\in \Ascr^\xi_{0}-\{0\}$;
\item $\alpha^i <\gamma <\alpha$ \, and \, $\Mod{\gamma}\leq \maxord{\xi}$, \,
for every $\gamma\in \Ascr^\xi_{\alpha^i}-\{\alpha^i\}$; and 
\end{enumerate}
for any given  maps $f\in Y^X$ and $g\in X^Z$, and any
multi-vector  $u\in Z^k$,
\begin{equation}
\label{tda}
 T_\alpha f \, \overline{u} = \sum_{\xi = \{\alpha^1,\ldots,\alpha^r\}\in {\mathscr P}_\alpha}
\Delta^{\,r}_{ u^\xi_{\alpha^1},\,\ldots,\, u^\xi_{\alpha^r} } f(u^{\xi}_0) \;,
\end{equation}
where
$u^{\xi}_\beta = \sum_{\gamma\in \Ascr^\xi_\beta} u_\gamma$,\,
for each \, $\beta= 0, \alpha^1,\ldots, \alpha^r$.
\end{teor}

\dem
This proposition is proved by induction in $k$.
It is obvious when $\alpha\in I^k$ with $k\leq 2$.
Assume it holds when $k\leq n$.
Any given multi-index $\alpha\in I^{k+1}$ is of the
form $\alpha \Diamond 0$, or $\alpha \Diamond 1$.
Since the first case follows by induction hypothesis,
we now restrict our attention to the second case.
Given $\overline{w}\in T^{k+1}(X)$,  write it as a
pair $\overline{w}=\pair{u}{v}$, of $k$-cuboids $\overline{u},\,\overline{v}\in T^k(X)$.
Using lemma~\ref{lemsplit} we deduce
\begin{align*}
T_{\alpha\Diamond 1}(f)(\overline{w}) &= T_{\alpha\Diamond 1}(f)\pair{u}{v} 
= T_\alpha(f)(\overline{u}+\overline{v}) -T_\alpha(f) \overline{u}\\
&= \sum_{\xi = \{\alpha^1,\ldots,\alpha^n\}\in {\mathscr P}_\alpha} 
\Delta^{\,n}_{ u^\xi_{\alpha^1}+v^\xi_{\alpha^1},\,\ldots,\, u^\xi_{\alpha^n}+v^\xi_{\alpha^n} } f(u^{\xi}_0+v^{\xi}_0) 
\;-\;
\Delta^{\,n}_{ u^\xi_{\alpha^1},\,\ldots,\, u^\xi_{\alpha^n} } f(u^{\xi}_0)\\
&=   \sum_{\xi = \{\alpha^1,\ldots,\alpha^n\}\in {\mathscr P}_\alpha}  
\Delta^{n+1}_{v_0^\xi,u_{\alpha^1}^\xi,\,\ldots,\,u_{\alpha^n}^\xi} f(u_0^\xi) \; +\\
  & \; +\; \sum_{i=1}^n \Delta^n_{u_{\alpha^1}^\xi,\,\ldots,\,
u_{\alpha^{i-1}}^\xi,\,v_{\alpha^i}^\xi,\,(u_{\alpha^{i+1}}^\xi+v_{\alpha^{i+1}}^\xi),\,
\ldots,\,(u_{\alpha^n}^\xi+v_{\alpha^n}^\xi)} f(u_0^\xi+v_0^\xi+u_{\alpha^i}^\xi)\, 
\end{align*}
To finish the proof we just need to establish a one-to-one correspondence between
summands above and partitions in ${\mathscr P}_{\alpha\Diamond 1}$.
Using the notation introduced in lemma~\ref{diamond},
the partition $\tilde{\xi}_0$ is associated with the first summand,
while the partitions $\tilde{\xi}_i$ ($1\leq i\leq n$) are associated each with one of the
subsequent $n$ summands.
Making these identifications we arrive at the equation between the previous sum
and the following one:
$$ \sum_{\xi = \{\alpha^1,\ldots,\alpha^n\}\in {\mathscr P}_\alpha} \left(\,
\Delta^{n+1}_{w_{0\Diamond 1}^{\tilde{\xi}_0}, 
w_{\alpha^1\Diamond 0}^{\tilde{\xi}_0},\,\ldots,\,w_{\alpha^n\Diamond 0}^{\tilde{\xi}_0} } 
f(w_{0\Diamond 0}^{\tilde{\xi}_0}) \; +\; 
\sum_{i=1}^n \Delta^n_{w_{\alpha^1\Diamond 0}^{\tilde{\xi}_{\alpha^i}},\,
\ldots,\, w_{\alpha^{i-1}\Diamond 0}^{\tilde{\xi}_{\alpha^i}},\, 
w_{\alpha^{i}\Diamond 1}^{\tilde{\xi}_{\alpha^i}}
\ldots,\, w_{\alpha^n\Diamond 1}^{\tilde{\xi}_{\alpha^i}}
} f(w_{0\Diamond 0}^{\tilde{\xi}_{\alpha^i}})
\,\right)
$$
which, by lemma~\ref{diamond}, is equal to
$$
\sum_{\tilde{\xi}= \{\beta^1,\ldots,\beta^r\} \in {\mathscr P}_{\alpha\Diamond 1}}  
\Delta^{\,r}_{ w^{\tilde{\xi}}_{\beta^1},\,\ldots,\, w^{\tilde{\xi}}_{\beta^r} } f(w^{\tilde{\xi}}_0) \;.  $$
We assume that $u^{\xi}_\beta = \sum_{\gamma\in \Ascr^\xi_\beta} u_\gamma$, and similar relations
hold for $v^{\xi}_\beta$, $w^{\tilde{\xi}}_{\beta\diamond 0}$, etc.
Then, matching terms, we arrive at the following recursive equations
on the sets $\Ascr^{\tilde{\xi}_i}_\beta$:

$$  \left\{\begin{array}{ccl}
\Ascr^{\tilde{\xi}_0}_{0\Diamond 0} &=& \Ascr_0^\xi \diamond 0 \\
\Ascr^{\tilde{\xi}_0}_{0\Diamond 1} &=& \Ascr_0^\xi \diamond 1 \\
\Ascr^{\tilde{\xi}_0}_{\alpha^i\Diamond 0} &=& \Ascr_{\alpha^i}^\xi \diamond 0
\end{array}
\right.$$
and
$$\left\{
\begin{array}{ll}
\Ascr^{\tilde{\xi}_i}_{\alpha^i\Diamond 1} = \Ascr_{\alpha^i}^\xi \Diamond 1  \quad 
&\text{ if } j=i\\
\Ascr^{\tilde{\xi}_i}_{\alpha^j\Diamond 0} = \Ascr_{\alpha^j}^\xi\Diamond 0  \quad &\text{ if } j<i\\
\Ascr^{\tilde{\xi}_i}_{\alpha^j\Diamond 0} = \Ascr_{\alpha^j}^\xi\Diamond 0  \,\cup\,
 \Ascr_{\alpha^j}^\xi\Diamond 1 \quad  &\text{ if } j>i\\
\Ascr^{\tilde{\xi}_i}_{0\Diamond 0} = \Ascr_0^\xi\Diamond 0 \,\cup\, \Ascr_0^\xi\Diamond 1 \,\cup\, 
\Ascr_{\alpha^i}^\xi\Diamond 0  
\end{array}
\right. $$
where $\Ascr \diamond i = \{ \, \beta\diamond i \,:\, \beta\in \Ascr \,\}$ for $i=0,1$.
These recursive relations ensure that $T_{\alpha\Diamond 1}(f)(\overline{w})$
has the correct development~(\ref{tda}). 
It is now easy to check inductively that these sets satisfy 
the conditions 1., 2., 3. and 4.
\cqd

\bigskip

Theorem A in the introduction is a corollary of theorem~\ref{TeorTdiscrfexpr}.
Its proof is completly analogous to that of proposition~\ref{invder}
and, therefore, will be omitted.

\bigskip

An algorithm that produces explicit expressions for $T_\alpha f\overline{u}$
has been devised, which was also used to confirm the correction of the above
recursive definitions. Its source code can be retrieved from 
{\rm http://ptmat.ptmat.fc.ul.pt/$\sim$pduarte/Research/FiniteDif\-ferences/index.html}. 
The following formulas were computer generated by this package.

\begin{align*}
 T_{({1, 1})}f(\overline{u}) = \;& \Delta_{u_{1,2}} f\left( u_{0}+u_{2}+u_{1}\right) + \\ 
& \Delta^2_{u_{1},u_{2}} f\left( u_{0}\right)\\ 
\end{align*}

\begin{align*}
 T_{({1, 1, 1})}f(\overline{u}) = \;& \Delta_{u_{1,2,3}} f\left( u_{0}+u_{3}+u_{2}+u_{2,3}+u_{1}+u_{1,3}+u_{1,2}\right) + \\ 
& \Delta^2_{u_{1},u_{2,3}} f\left( u_{0}+u_{3}+u_{2}\right) + \\ 
& \Delta^2_{u_{1,3},u_{2}+u_{2,3}} f\left( u_{0}+u_{3}+u_{1}\right) + \\ 
& \Delta^2_{u_{1,2},u_{3}+u_{2,3}+u_{1,3}} f\left( u_{0}+u_{2}+u_{1}\right) + \\ 
& \Delta^3_{u_{1},u_{2},u_{3}} f\left( u_{0}\right)\\ 
\end{align*}

\begin{align*}
 \Delta^2_{v} (f \circ g )(x) = \;& \Delta_{\Delta^2_{v_1,v_2} g\left( x\right)} f\left( g(x)+\Delta_{v_1} g\left( x\right)+\Delta_{v_2} g\left( x\right)\right) + \\ 
& \Delta^2_{\Delta_{v_1} g\left( x\right),\Delta_{v_2} g\left( x\right)} f\left( g(x)\right)\\ 
\end{align*}

\begin{align*}
 \Delta^3_{v} (f \circ g )(x) = \;& \Delta_{\Delta^3_{v_1,v_2,v_3} g\left( x\right)} f\left( g(x)+\Delta_{v_1} g\left( x\right)+\Delta_{v_2} g\left( x\right) + \Delta_{v_3} g\left( x\right)\;+ \right. \\
& \qquad \qquad \qquad \qquad \left. 
\Delta^2_{v_1,v_2} g\left( x\right)+\Delta^2_{v_1,v_3} g\left( x\right)+\Delta^2_{v_2,v_3} g\left( x\right)\right) \; + \\ 
& \Delta^2_{\Delta_{v_1} g\left( x\right),\Delta^2_{v_2,v_3} g\left( x\right)} f\left( g(x)+\Delta_{v_2} g\left( x\right)+\Delta_{v_3} g\left( x\right)\right) \; + \\ 
& \Delta^2_{\Delta^2_{v_1,v_3} g\left( x\right),\Delta_{v_2} g\left( x\right)+\Delta^2_{v_2,v_3} g\left( x\right)} f\left( g(x)+\Delta_{v_1} g\left( x\right)+\Delta_{v_3} g\left( x\right)\right) \; + \\ 
& \Delta^2_{\Delta^2_{v_1,v_2} g\left( x\right),\Delta_{v_3} g\left( x\right)+\Delta^2_{v_1,v_3} g\left( x\right)+\Delta^2_{v_2,v_3} g\left( x\right)} f\left( g(x)+\Delta_{v_1} g\left( x\right)+\Delta_{v_2} g\left( x\right)\right) \; + \\ 
& \Delta^3_{\Delta_{v_1} g\left( x\right),\Delta_{v_2} g\left( x\right),\Delta_{v_3} g\left( x\right)} f\left( g(x)\right)
\end{align*}

\bigskip

\thispagestyle{empty}

\end{document}